\documentclass[12pt]{article}%
\usepackage{amssymb}
\usepackage{amsfonts}
\usepackage[utf8]{inputenc}
\usepackage{amsmath}
\usepackage{graphicx}%
\usepackage { hyperref}
\providecommand{\U}[1]{\protect\rule{.1in}{.1in}}
\newtheorem{theorem}{Theorem}[section]

\newtheorem{defin}[theorem]{Definition}

\newtheorem{lem}[theorem]{Lemma}
\newtheorem{teo}[theorem]{Theorem}

\newtheorem{propo}[theorem]{Proposition}
\newtheorem{remark}[theorem]{Remark}

\usepackage{authblk}

\begin{document}

\title{The $p-$parabolicity under a decay assumption on the Ricci curvature}
\author{L. S. Priebe, R. B. Soares}
\date{}
\maketitle

\begin{abstract}
We prove that, given $\alpha>0$, if $M$ is a complete Riemannian manifold which Ricci
curvature satisfies.
\begin{equation}
\operatorname*{Ric}\nolimits_{x}(v)\geq\alpha\operatorname{sech}^{2}%
(r(x)))\label{8}%
\end{equation}
or
\begin{equation}
\operatorname*{Ric}\nolimits_{x}(v)\geq-\frac{{h_{\alpha}}%
(r(x))}{r(x)^{2}},\label{9}%
\end{equation}
where
\[
{h_{\alpha}}(r)=\frac{\alpha(\alpha+1)r(x)^{\alpha
}}{r(x)^{\alpha
}-1},
\]
for all $x\in M\backslash B_{R}(o)$ and for all $v\in T_{x}M,$
$\left\Vert v\right\Vert =1,$ where \ $o$ is a fixed point of $M$,
$r(x)=d(o,x)$, $d$ the Riemannian distance in $M$ and $B_{R}(o)$ the geodesic
ball of $M$ centered at $o$ with radius $R>0$, then $M$ is $p-$parabolic for any
$p>1$, if satisfies $(\ref{8})$, and $M$ is $p-$parabolic, for any
$p\geq(\alpha+1)(n-1)+1$ if satisfies $(\ref{9})$.

\end{abstract}

\section{Introduction}

\qquad Let $M$ be a complete Riemannian manifold. Given $p>1,$ the
$p-$Laplacian $\Delta_{p}$ in $M$ is defined by
\[
\Delta_{p}(u)=\operatorname{div}\left(  \frac{\nabla u}{\left\Vert \nabla
u\right\Vert ^{2-p}}\right)  ,\text{ }u\in C^{2}\left(  M\right)  .
\]
For $p=2,$ $\Delta=\Delta_{2}$ is the usual Laplacian operator.

A function $u\in C^{1}\left(  M\right)  $ is called $p-$superharmonic if
$\Delta_{p}(u)\leq0$ (in the weak sense). $M$ is called $p-$parabolic if any
positive $p-$superharmonic (entire) function is constant.

Yau proved (Corollary 3.1 of \cite{Y}) that if the Ricci curvature of $M$ is
nonnegative then $M$ satisfies the Liouville property that is, any positive
harmonic function in $M$ is constant (see \cite{K}). Liouville property does
not extend to $2-$parabolicity that is, one cannot replace harmonic by
superharmonic in Yau's result since $\mathbb{R}^{n},$ $n\geq3,$ is well known
to be non $2-$parabolic. Nevertheless, if the Ricci curvature, besides being
non negative, does not converge too fast to zero, this extension holds,
including the $p-$parabolicity. 

Choose a point $o\in M$ and set $r(x):=d(o,x),$ $x\in M,$ where $d$ is the
Riemannian distance in $M.$ Precisely, we prove:

\begin{teo}
\label{th}Let $M$ be a complete Riemannian manifold. Given $\alpha>0$, assume
that the Ricci curvature of $M$ satisfies
\[
\operatorname*{Ric}\nolimits_{x}(v)\geq\alpha\operatorname{sech}^{2}(r(x)),
\]
for all $x\in M\backslash B_{R}(o)$, all $v\in T_{x}M,$ $\left\Vert
v\right\Vert =1,$ and for some $R>0.$ Then $M$ is $p-$parabolic, for any
$p>1.$
\end{teo}

J. B. Casteras, E. Heinonen and I. Holopainen, in \cite{WH}, prove that, given $\alpha>0$, if%
\[
{K}(x)\geq-\frac{\alpha}{r(x)^{2}\operatorname*{log}(r(x))},
\]
for any $x$ outside a compact subset of $M$ where, for a given $x,$
${K}(x)$ is the infimum of the sectional curvatures of $M$ on
the planes through $T_{x}M,$ then $M$ is $p$-parabolic if $p=n$ and
$0<\alpha\leq1$ or $p>n$ and $\alpha>0$. In \cite{A} A. V\"{a}h\"{a}kangas
proves that if%
\[
{K}(x)\geq-\frac{\alpha(\alpha+1)}{r(x)^{2}},
\]
then $M$ is $p$-parabolic if $p\geq(\alpha+1)(n-1)+1$. 

To see more clearly the relation between our next theorem and Vahakangas' result it is conveneint to introduce the function $h:(1,\infty)\rightarrow \mathbb{R}$,
\[
{h_{\alpha}}({r)}=\frac{\alpha(\alpha+1)r^{\alpha}%
}{r^{\alpha}-1%
}=\alpha(\alpha+1)\frac{r^{\alpha}%
}{r^{\alpha}-1%
}>\alpha(\alpha+1).
\]
Hence, the next result improves V\"{a}h\"{a}kangas' on the condition on the
decay of the curvature and also by replacing the sectional curvature by the
Ricci curvature:

\begin{teo}
\label{th2}Let $M$ be a complete Riemannian manifold. Given $\alpha>0$, define
\begin{equation}\label{12}
    {h_{\alpha}}(r)=\frac{\alpha(\alpha+1)r^{\alpha}}{r^{\alpha}-1}.
\end{equation}

Assume that the Ricci curvature of $M$ satisfies
\[
\operatorname*{Ric}\nolimits_{x}(v)\geq-\frac{{h_{\alpha}(r(x))}%
}{r(x)^{2}},
\]
for all $x\in M\backslash B_{R}(o)$, all $v\in T_{x}M,$ $\left\Vert
v\right\Vert =1,$ and for some $R>0.$ Then $M$ is $p-$parabolic, for any
$p\geq(\alpha+1)(n-1)+1$.
\end{teo}

It is interesting to mention that, closely related to Theorem \ref{th2}, there
is a conjecture of Green-Wu \cite{L1} which goes in an opposite direction,
namely: if $\operatorname*{K}\leq-C/r^{2}$ outside a compact subset of $M,$
then there are non constant, bounded, harmonic functions on $M.$ Here
$\operatorname*{K}(x),$ $x\in M,$ is the supremum of the sectional curvatures
of $M$ on the planes through $T_{x}M.$

The authors are indebted with Professor Detang Zhou who suggested Lemmas
\ref{l}, \ref{l2}, and their proofs.

\section{Preliminaries}

\qquad For proving Theorem \ref{th} we use a well known criterion obtained by Ilkka
Holopainen, namely:

\begin{theorem}
\label{w} \label{il}Let $M$ be a complete Riemannian manifold. Fix a point
$o\in M$ and write $V(r)=Vol(B(o,r))$. Given $1<p<\infty$. If
\[
\int^{\infty}\left(  \frac{r}{V(r)}\right)  ^{\frac{1}{p-1}}dr=\infty
\]
or
\[
\int^{\infty}\left(  \frac{1}{V^{\prime}(r)}\right)  ^{\frac{1}{p-1}}%
dr=\infty,
\]
then $M$ is $p-$parabolic.
\end{theorem}

To make use of Theorem \ref{il} we introduce a definition and prove two lemmas.

\begin{defin}
(Riccati equation and solution) Let $M$ be a Riemannian manifold. Let
$\gamma:I\rightarrow M$ be a unit speed geodesic and $U\in End(\gamma^{\perp
})$ be a symmetric field of endomorphisms along $\gamma$. If ${R}:={R}_{v}%
:={R}(\text{ }_{,}\gamma^{\prime})\gamma^{\prime\perp},\text{ }\in\gamma^{\perp
}$, $v=\gamma^{\prime}(0)$, is the curvature endomorphism along $\gamma$, we
call
\[
{U} ^{\prime}+{U}^{2}+{R}=0
\]
the Riccati equation. We say ${U}\in End(\gamma^{\perp})$ is a solution of the
Riccati equation, if
\[
\forall X\in\mathcal{X}(\gamma^{\perp}):{U}^{\prime}(X)+{U}^{2}(X)+{R}(X)=0.
\]
\end{defin} 

\begin{theorem}\label{t2.3}
    (Riccati and Jacobi equation). Let $c:I \rightarrow M$ be a unit speed geodesic, $E_{1} = c^{\prime}$, and
let $E_{2}, \ldots , E_{n}$ be a parallel ONB along $c$. Let $t_{0} \in I$ and let $X_{2}, \ldots , X_{n}$ be any basis of $(c^{\prime}(t_{0}))^{\perp}$. For 
any $2 \leq i \leq n$ let $J_{i}$ be the Jacobi field along $c$ satisfying
\[
J_{i}(t_{0}) = X_{i}\  \text{and} \     D_{t}J_{i}(t_{0}) = U_{0}(X_{i}),
\]
where $U_{0}$ is a given symmetric endomorphism on $c^{\prime}(t_{0})^{\perp}$. Define a tensor $J_{t}: c^{\prime}(t_{0})^{\perp} \rightarrow c^{\prime}(t_{0})^{\perp}$ along $c$ by 
\begin{equation}\label{eq22}
    J_{t}\left(\sum_{i=2}^{n} \alpha_{i} E_{i} \right) :=\sum_{i=2}^{n} \alpha_{i} J_{i}.
\end{equation}

The endomorphism $J = J_{t}$ solves the Jacobi equation $J^{\prime\prime}+RJ=0$, is invertible for any $t$ near $t_{0}$ and
\[
U_{t}:=J_{t}^{\prime}\circ J_{t}^{-1}.
\]
is a symmetric solution of the Riccati equation.
\end{theorem}

\textbf{Proof:} See \cite{B}.

\bigskip

\begin{lem}
\label{l} Let $M$ be a Riemannian manifold. Let $\gamma:I\rightarrow M$ be a
unit speed geodesic, $U$ be a symmetric solution of the Riccati equation and
$J$ be a field of isomorphisms satisfying $J^{\prime}=J\circ U$ U (in particular one may choose $U$ and
$J$ as in Theorem \ref{t2.3}). Define
$u:I\rightarrow M$ by
\[
u:=\frac{1}{n-1}\frac{\partial}{\partial t}(\ln(\det(J))).
\]
Then
\[
u=\frac{1}{n-1}\operatorname*{tr}(J)
\]
and
\[
u^{\prime}\leq-u^{2} -\frac{1}{n-1}\operatorname*{Ric}(\gamma^{\prime},\gamma^{\prime}).
\]

\end{lem}

\textbf{Proof:} See \cite{B}.
\bigskip

Fixed a point $o\in M$, in terms of the polar normal coordinates at $o$, we
can write the volume element as $J^{n-1}drd\theta$.

\begin{lem}
\label{l2} Let $(M^{n},g)$ be a complete Riemannian manifold with Ricci
curvature bounded from below by $\frac{(n-1)}{2}{K}$. Then, for any $x=(\theta,r)$ that is
not in the cut-locus $C(p)$ of $p$, we have
\[
\ln\left(  \frac{J}{\phi}\right)  \leq-\int_{{0}}^{r}\left(  \frac{1}{\phi
(t)}\int_{{0}}^{t}(\phi^{\prime\prime}(s)-\phi(s){K}(s))ds\right)  dt,
\]
where $J^{n-1}$ is the volume element and $\phi:[0,l)\longrightarrow
\mathbb{R}$ is a function satisfying $\phi(0)$, $\phi^{\prime}(0)=1$ and
$\phi(r)>0$, for $r\in(0,l)$.
\end{lem}

\textbf{Proof: }
Firstly, we observe that
\begin{eqnarray*}
    Vol(B_{o}(R))=\int_{\mathbb{S}^{n-1}} \int_{0}^{\min{R,cut{\theta}}} J^{n-1} drd\theta.
\end{eqnarray*}

Considering $\operatorname*{exp_{o}}(t,\theta): [0,\infty) \rightarrow M$ the exponential map, follow that the coordinated fields are $X_{1}=\frac{\partial}{\partial r }\operatorname*{exp_{o}}(t,\theta)$ and $X_{i}=\frac{\partial}{\partial \theta_{i}}\operatorname*{exp_{o}}(t,\theta)$, for $i>1$, and  $g_{ij}$ is the metric related to this parameterization.

Note that $g_{11}=1$ because it is geodesic and by Lemm Gauss $g_{1i}=0$, for $i>1$, therefore $ J^{n-1}=\det(g_{i,j})$, for $i,j>1$.

Given $q\in C(o)$, the set C(o) of all cut points of o is
the cut locus of o, there is $\theta$ such that $\gamma(t)=\operatorname*{exp_{o}}(t,\theta)$ it is the geodesic that $\gamma(0)=o$ and  $|\gamma'(0)|=1$. 

Note that $J_{i}(t)=X_{i}(t,\theta)$ is Jacobi field with $J(0)=0$, $J'(0)=1$ along the $\gamma$. Furthermore, by the fact that $q\in C(o)$, follow that
$J_{i}\perp J_{j}$. Then, defined $E_{i}$, for $i>1$, and $J_{t}$ according to Teorem \ref{t2.3}, we have
\begin{eqnarray*}
   det(J_{t})=det(g(J_{t}(E_{i}),{J_{j}}))=det(g(J_{i},{J_{j}}))=\det(g_{i,j})_{n\times n}.
\end{eqnarray*}

Therefore, writing $J(t)=J(t,\theta)^{n-1}$, we have that%
\[
w=\frac{\partial}{\partial r}\ln((\det{J})),
\]
satisfies the hypotheses of Lemma \ref{l}. Then,
\[
w^{\prime}+w^{2}+{K}\leq0.
\]

Multiplying the above inequality by $(\phi(t))^{2}$ and integrating from $0$ to
$r$, we have to
\[
\int_{0}^{r}(\phi(t)^{2}w^{\prime}(t)+\phi(t)^{2}w(t)^{2}+\phi(t)^{2}{K}(t))dt\leq0.
\]

Using $\phi({0})=0$,
\begin{align*}
\phi(r)^{2}w(r)  &  =\int_{{0}}^{r}\frac{\partial}{\partial t}(\phi
(t)^{2}w(t))dt\\
&  =\int_{{0}}^{r}2\phi(t)\phi^{\prime}(t)w(t)dt+\int_{{0}}^{r}\phi
(t)^{2}w^{\prime}(t)dt.
\end{align*}

Then
\[
\int_{{0}}^{r}\phi
(t)^{2}w^{\prime}(t)dt=\phi(r)^{2}w(r)-\int_{0}^{s}2\phi(t)\phi^{\prime
}(t)w(t)dt
\]
and consequently
\[
\phi(r)^{2}w(r)+\int_{{0}}^{s}(-2\phi(t)\phi^{\prime}(t)w(t)%
+\phi(t)^{2}w(t)^{2}+\phi(t)^{2}{K}(t))dt\leq0.
\]

In this way,
\begin{align*}
\phi(r)^{2}w(r) &  \leq-\int_{{0}}^{r}(-2\phi(t)\phi^{\prime}(t)%
w(t)+\phi(t)^{2}w(t)^{2}+\phi(t)^{2}{K}(t))dt\\
&  \leq-\int_{{0}}^{r}(-2\phi(t)\phi^{\prime}(t)w(t)+\phi(t)^{2}w(t)^{2}+\phi^{\prime}(t)^{2}%
-\phi^{\prime}(t)^{2}+\phi(t)^{2}{K}(t))dt\\
&  \leq-\int_{{0}}^{r}(\phi(t)w(t)-\phi^{\prime}(t))^{2}dt  +\int_{{0}}^{r}(\phi^{\prime}(t)^{2}-\phi(t)^{2}{K}(t))dt\\
&  \leq \int_{{0}}^{r}(\phi^{\prime}(t)^{2}-\phi(t)^{2}{K}(t))dt.
\end{align*}

Observe that
\begin{align*}
\phi(r)\phi^{\prime}(r) &  =\int_{{0}}^{r}\frac{\partial}{\partial t}%
(\phi(t)\phi^{\prime}(t))dt\\
&  =\int_{{0}}^{r}\phi^{\prime}(t)^{2}dt+\int_{{0}}^{r}\phi(t)\phi^{\prime\prime
}(t)dt,
\end{align*}
that is,
\[
\int_{{0}}^{r}\phi^{\prime}(t)^{2}dt=\phi(r)\phi^{\prime}(r)-\int_{{0}}^{r}%
\phi(t)\phi^{\prime\prime}(t)dt.
\]
It follows that
\[
\phi(r)^{2}w(r)\leq\phi(r)\phi^{\prime}(r)+\int_{{0}}^{r}(-\phi
(t)\phi^{\prime\prime}-\phi(t)^{2}{K}(t))dt.
\]
Then
\[
w(r)\leq\frac{\phi^{\prime}(r)}{\phi(r)}-\frac{1}{\phi(r)^{2}}\int_{{0}}^{r}%
\frac{1}{\phi(t)}(\phi^{\prime\prime}(t)+\phi(t){K}(t))dt.
\]

Moreover, observe that
\[
\lim_{r\longrightarrow0}\frac{J(r)}{\phi(r)}=\lim_{r\longrightarrow0}%
\frac{J^{\prime}(r)}{\phi^{\prime}(r)}=\frac{1}{1}=1.
\]
Therefore
\[
\ln\left(  \frac{J}{\phi}\right)  (0)=1.
\]

Then
\begin{align*}
 \ln\left(  \frac{J}{\varphi}\right)  (r) & =\int_{{0}}^{r}\frac{\partial}{\partial t}\left(  \ln\left(  \frac
{J(t)}{\phi(t)}\right)  \right)  dt\\
&  =\int_{{0}}^{r}\frac{\partial}{\partial t}\left(  \ln\left(  {J(t)}\right)
-\ln\left(  {\phi(t)}\right)  \right)  dt\\
&  \leq\int_{{0}}^{r}\left(  \frac{\phi^{\prime}(t)}{\phi(t)}-\frac{1}%
{\phi(t)}\int_{{0}}^{t}(\phi^{\prime\prime}(s)-\phi(s)K(s))ds-\frac
{\phi^{\prime}(t)}{\phi(t)}\right)  dt\\
&  \leq-\int_{{0}}^{r}\left(  \frac{1}{\phi(t)}\int_{{0}}^{t}(\phi
^{\prime\prime}(s)-\phi(s){K}(s))ds\right)  dt.
\end{align*}

What concludes the proof of the lemma. We may now prove the following theorem:

\begin{teo}
\label{a1} Let $(M^{n},g)$ be a complete Riemannian manifold and
${K}:[0,\infty)\longrightarrow\mathbb{R}$ be a continuous function. If there
is $o\in M$ such that
\[
\operatorname*{Ric}\nolimits_{x}(v)\geq(n-1){K}(r(x)),
\]
for all $x\in M$, $v\in T_{x}M$, where $r(x)=d(x,o)$, then, the volume of the geodesic ball
$B_{p}(r)$ satisfies
\[
\operatorname*{Vol}(B_{p}(r))\leq w_{n-1}\int_{0}^{r}\phi(r)^{n-1}dr,
\]
where $w_{n-1}$ is the volume of the $(n-1)-$dimensional sphere end
$\phi:[0,\infty)\longrightarrow\mathbb{R}$ is the solution of
\begin{align*}
\phi^{\prime\prime}(t)+{K}(t)\phi(t)  =0\\
\phi(0)  =0\ \phi^{\prime}(0)=1.
\end{align*}

\end{teo}

\textbf{Proof: }

Note that,
\[
Vol(B_{p}(r))=\int_{\theta}\int_{0}^{r}J^{n-1}(r,\theta)drd\theta.
\]

Moreover, by Lemma \ref{l2}
\begin{align*}
\ln\left(  \frac{J}{\phi}\right)   &  \leq-\int_{0}^{r}\frac{1}{\phi^{2}%
(t)}\int_{0}^{s}(\phi(t)\phi^{\prime\prime2}\phi^{2}(t))dsdt\\
&  \leq0.
\end{align*}

Therefore
\begin{align*}
\frac{J}{\phi}\leq1,
\end{align*}
consequently,
\begin{align*}
{J}\leq{\phi}.
\end{align*}

Then
\begin{align*}
Vol(B_{p}(r)) &  =\int_{\theta}\int_{0}^{r}J^{n-1}(r,\theta)drd\theta\\
&  \leq\int_{\theta}\int_{0}^{r}\phi^{n-1}(r)drd\theta\\
&  =w_{n-1}\int_{0}^{r}\phi^{n-1}(r)drd\theta.
\end{align*}

Finally, the following Lemma results from a real analysis. It can be
demonstrated using the techniques described in Lemma $2.31$ in \cite{L1} and Lemma
$5$ in \cite{M}:

\begin{propo}
\label{p1} Let $f,g\in C^{2}([R,\infty))$, $R\geq0$, such that
\begin{align*}
-\frac{g^{\prime\prime}}{g}= -\frac{f^{\prime\prime}}{f},
\end{align*}
where $\alpha> 0$. If there is $R_{0}>0$ such that $f(R_{0})\geq0$ and
$f^{\prime}(R_{0})\geq0$, then there is $\beta>0$ such that
\begin{align*}
g(r) \leq\beta f(r),
\end{align*}
for $r>R_{0}$.
\end{propo}

\section{Proof of Theorem \ref{th}}

Firstly, for $\alpha>0$, take $\delta>0$ such that $\delta<
\frac{\alpha}{2}$, defining
\begin{align*}
K_{\delta}(r)=-\delta\operatorname{csch}^{2}(r) (-1 + \delta
\operatorname{sech}^{2}(r) -\operatorname{tanh}^{2}(r)),
\end{align*}
we have to
\begin{align*}
\lim_{r\rightarrow\infty}\frac{K_{\delta}(r)}{\alpha\operatorname{sech}%
^{2}(r)}=\frac{2 \delta}{\alpha}<1.
\end{align*}
Then, there is $r_{\delta} > R$ such that for $r > r_{\delta}$
we have
\begin{align*}
K_{\delta}(r)\leq\operatorname{sech}^{2}(r),
\end{align*}
and consequently
\begin{align*}
\operatorname{Ric}_{x}(v)\geq\alpha K_{\delta}(r(x)),
\end{align*}
for $x\in M-B_{r_{o}}(r_{\delta})$.

Let
\[
B=inf\{Ric_{x}(v);\forall x\in M-B_{p}(r_{\delta}),\forall v\in
T_{x}M\},
\]
$\overline{K}:\mathbb{R}_{+}\longrightarrow\mathbb{R}$ a smooth function
defined by
\[
\overline{K}=\left\{
\begin{array}
[c]{lll}%
B,\ \ if\ \ r<r_{\delta} &  & \\
f(r),\ \ if\ \ r_{\delta}\leq r\leq r_{\delta}+1 &  & \\
K_{\delta}(r),\ \ if\ \ r_{\delta}+1\leq r &  &
\end{array}
\right.  ,
\]
where $f(r)$ is a smooth function with $f(r)\leq K_{\delta}(r)$, for
$r\in\lbrack r_{\delta},r_{\delta}+1]$.

By Existence and Uniqueness Theorem, there is $\overline{\phi}(r)$ such that
\begin{align*}
\overline{\phi}^{\prime\prime}(r)+\overline{K}\overline{\phi}(r)=0\\
\overline{\phi}(0)=0\ \overline{\phi}^{\prime}(0)=1.
\end{align*}

Using the Theorem \ref{a1}, $Vol(B_{p}(r))$ the volume of the geodesic ball,
follow that
\begin{align*}
Vol(B_{p}(r))  &  \leq w_{n-1}\int_{0}^{r}(\overline{\phi}(t))^{n-1}dt.
\end{align*}

Note that a solution of ode
\begin{align*}
{\phi}^{\prime\prime}(r)+K_{\delta}{\phi}(r)=0
\end{align*}
is given by $\phi(r)=\operatorname{tanh}^{\delta}(r)$.

By Proposition \ref{p1} there is $\beta^{\frac{1}{n-1}}>0$ such that
$\overline{\phi}(r)\leq\beta^{\frac{1}{n-1}}\operatorname{tanh}^{\delta
}(r)$, for $r>r_{\delta}+1$. Then, $r>r_{\delta}+1$,
\begin{align*}
Vol(B_{p}(r))  &  \leq w_{n-1}\int_{0}^{r}(\overline{\phi}(t))^{n-1}dt\\
&  \leq w_{n-1}\left(  \int_{0}^{r_{\delta}+1}(\overline{\phi}%
(t))^{n-1}dt+\int_{r_{\delta}+1}^{r}(\beta^{\frac{1}{n-1}}%
\operatorname{tanh}^{\delta}(r))^{n-1}dt\right) \\
&  \leq w_{n-1}\left(  \int_{0}^{r_{\delta}+1}(\overline{\phi}%
(t))^{n-1}dt+\int_{r_{\delta}+1}^{r} \beta dt\right) \\
&  \leq w_{n-1}\left(  \int_{0}^{r_{\delta}+1}(\overline{\phi}%
(t))^{n-1}dt+\beta({r} -(r_{\delta}+1)) \right) \\
&  \leq C r,
\end{align*}
where $C$ is a positive constant.

This way
\[
\int^{\infty}\left(  \frac{t}{Vol(B({t}))}\right)  ^{\frac{1}{p-1}}dt\geq
\int^{\infty}\left(  \frac{t}{Ct}\right)  ^{\frac{1}{p-1}}dt=\infty,
\]
for all $p>1$. Then, by Theorem $\ref{w}$, $M$ is\ $p-$parabolic for all $p>1$.

\section{Proof of Theorem \ref{th2}}

Equivalent to the proof of Theorem \ref{th}, we consider
\[
K(r)=h_{\alpha}(r)\frac{1}{r^{2}},
\]
and defining $\overline{K}:\mathbb{R}_{+}\longrightarrow\mathbb{R}$, a smooth
function, by
\[
\overline{K}=\left\{
\begin{array}
[c]{lll}%
B,\ \ if\ \ r<R &  & \\
f(r),\ \ if\ \ R\leq r\leq R+1 &  & \\
K(r),\ \ if\ \ R+1\leq r &  &
\end{array}
\right.  ,
\]
where
\[
B=inf\{Ric_{x}(v);\forall x\in M-B_{p}(R),\forall v\in T_{x}M\},
\]
and $f(r)$ is a smooth function with $f(r)\leq K(r)$, for $r\in\lbrack R,R+1]$.

Note that $\phi(r)=r^{\alpha+1}-r$ is solution of
\begin{align*}
{\phi}^{\prime\prime}(r)+K{\phi}(r)=0.
\end{align*}
Using the Theorem \ref{a1} and Proposition \ref{p1}, there is $
{\beta}^{\frac{1}{n-1}}>0$ such that for $r>R+1$
\begin{align*}
Vol(B_{p}(r))  &  \leq Vol(B_{p}(R+1))+w_{n-1}\int_{R+1}^{r}\left(
{\beta}^{\frac{1}{n-1}}(r^{\alpha+1}-1)\right)
^{n-1}dt\\
&  \leq Vol(B_{p}(R+1))+\beta w_{n-1}\int_{R+1}^{r}(r^{\alpha+1})^{n-1}dt\\
&  \leq C r^{({\alpha+1})(n-1)+1},
\end{align*}
where $C$ is a positive constant.

This way
\begin{align*}
\int^{\infty}\left(  \frac{t}{Vol(B({t}))}\right)  ^{\frac{1}{p-1}}dt\geq
\int^{\infty}\left(  \frac{t}{C t^{({\alpha+1})(n-1)+1}}\right)  ^{\frac
{1}{p-1}}dt=\infty,
\end{align*}
if $p\geq(\alpha+1)(n-1)+1$.

Then, by Theorem $\ref{w}$, $M$ is\ $p-$parabolic for all $p\geq
(\alpha+1)(n-1)+1$.

\begin{remark}
It is important to note that the function ${h_{\alpha}}$ in the $\operatorname{Ricci}$ curvature hypothesis of Theorem \ref{th2}, it admits an even more general definition, which will be indicatd by $h_{\alpha,\epsilon}$, with $\alpha$ and $\varepsilon$ positive constants and $\alpha+1>\epsilon$, which also guarantees the validity of the result, that is, if
\begin{equation}
\operatorname*{Ric}\nolimits_{x}(v)\geq-\frac{{h_{\alpha,\varepsilon}}%
(r(x))}{r^{2}(x)},\label{9}%
\end{equation}
where
\[
{h_{\alpha,\varepsilon}}(r)=\frac{\alpha(\alpha+1)(r(x))^{\alpha
}-\varepsilon(\varepsilon-1)(r(x))^{\varepsilon
}}{(r(x))^{\alpha
}-(r(x))^{\varepsilon
}},
\]
for all $x\in M\backslash B_{R}(o)$, all $v\in T_{x}M,$ $\left\Vert
v\right\Vert =1,$ and for some $R>0$, then $M$ is $p-$parabolic, for any
$p\geq(\alpha+1)(n-1)+1$.
\end{remark}


\begin{thebibliography}{9}                                                                                                %


\bibitem {H}I. Holopainen, Quasiregular mappings and the p-Laplace operator,
Contemp. Math. 338 (2003), 219-239.

\bibitem {K}J. Kazdan: Parabolicity and the Liouville Property on Complete
Riemannian Manifolds, eminar on New Results in Nonlinear Partial Differential
Equations, Aspects of Mathematics, V. 10, 153-166, 1987

\bibitem {P}P. Petersen, Riemannian Geometry, Graduate Texts in Mathematics,
Springer-Verlag, 2006.

\bibitem {B}W. Ballmann, Differential Geometry II, 2010. Avaiable at:\url{https://luis.impa.br/aulas/georiem/Ballman_DifferentialGeometryII.pdf},
access in 10/10/23.

\bibitem {Y}R. Schoen and S. Yau. Lectures on Harmonic Maps. Conference
Proceedings and Lecture Notes in Geometry an Topology, v. II, p. 394, 1997.


\bibitem {Y}R. Schoen and S. Yau. Lectures on Harmonic Maps. Conference
Proceedings and Lecture Notes in Geometry an Topology, v. II, p. 394, 1997.

\bibitem {M}P. March.Brownian Motion and Harmonic Functions on Rotationally
Symmetric Manifolds.The Annals of Probability, v. 14, pp. 793-801, 1986.

\bibitem {L1}R. E. Greene , H. Wu. Lecture Notes in Mathematics.Function
Theory on Manifolds Which Possess a Pole, v. 699, 1979.

\bibitem {WH}J.B. Casteras, E. Heinonen, I. Holopainen. Existence and
non-existence of minimal graphic and p-harmonic functions, Proc. Roy. Soc.
Edinburgh Sect. A 150 (2020), 341-366

\bibitem {A}V\"{a}h\"{a}kangas, A. Dirichlet Problem at Infinity for
A-Harmonic Functions. Potential Anal 27, 27 - 44 (2007). https://doi.org/10.1007/s11118-007-9051-7
\end{thebibliography}
\end{document}